\documentclass[61pt]{article}

\usepackage{amssymb}
\usepackage{graphicx}
\usepackage{amsmath}
\usepackage{amsfonts}
\usepackage{latexsym}
\usepackage{epstopdf}
\usepackage{float}
\usepackage{subfigure}

\allowdisplaybreaks[4]
\usepackage{amsmath}
\allowdisplaybreaks[4]
\newtheorem{theorem}{{\bf Theorem}}[section]
\newtheorem{lemma}[theorem]{{\bf Lemma}}
\newtheorem{prop}[theorem]{{\bf Proposition}}

\newenvironment{proof}{\noindent{\em Proof. }}{\quad \hfill $\Box$\vspace*{2ex}}
\newtheorem{definition}[theorem]{{\bf Definition}}

\newenvironment{proof th1.1}{\noindent{\em Proof\textbf{( of Theorem 1.1)}. }}{\quad \hfill $\Box$\vspace*{2ex}}
\newenvironment{proof th1.3}{\noindent{\em Proof\textbf{( of Theorem 1.3)}. }}{\quad \hfill $\Box$\vspace*{2ex}}
\newenvironment{proof th1.4}{\noindent{\em Proof\textbf{( of Theorem 1.4)}. }}{\quad \hfill $\Box$\vspace*{2ex}}

\begin{document}
\baselineskip 20pt

    \title{{  Dynamical Sampling in
Shift-Invariant Spaces Associated with multi-dimensional  Special Affine Fourier Transform }\footnote{Supported by NSFC (No.11401435) }}
   \author{Meng Ning, Li-Ping Wu\thanks{Corresponding author.
   E-mail: wuliping@email.tjut.edu.cn},  Qing-yue Zhang, Bei Liu
    \\
    \\
School of Science\\
Tianjin  University of Technology\\
Tianjin 300384\\
People's Republic of China\\
\\
}
\date{}
  \maketitle
  \begin{center}
\begin{minipage}{145mm}
{\small{\mbox{}\hspace{8pt}

\mbox{}\hspace{8pt}
 }

\mbox{}\hspace{8pt}}
 \end{minipage}
 \end{center}
 \begin{abstract}

  The Special Affine Fourier Transformation(SAFT), which generalizes several well-known unitary transformations, has been demonstrated as a valuable tool in signal processing and optics.  In this paper, we explore the multivariate dynamical sampling problem
in shift-invariant spaces associated with the multi-dimensional SAFT.  Specifically, we derive a sufficient and necessary
condition under which  a function in a shift-invariant space can be stably recovered from its dynamical sampling measurements
associated with the multi-dimensional SAFT . We also present
a  straightforward example to elucidate our main result.

 \textbf{2010 MR Subject Classification :} 94A20, 94A12, 42C15

 \end{abstract}

\bigskip
\section{Introduction}

\hskip \parindent
The one-dimensional Special Affine Fourier Transformation (SAFT), which was introduced in \cite{sa},  is an integral transformation
associated with a general inhomogeneous lossless linear mapping in phase-space. It has been shown to be powerful tool for signal processing and optics \cite{sc,scj}. In this paper, we define multi-dimensional SAFT, which includes many well-known linear  transforms as its special cases, such as the Fourier transform (FT),
the offset FT, the fractional Fourier transform (FRFT), the Fresnel transform,  the linear canonical transform(LCT), The Lorentz transform(LOT) described in Table I below \cite{zhang,sly,l}.
  Therefore, investigating relevant theorems for multi-dimensional SAFT may help to improve  more understanding on its special cases and to gain many applications across the fields of signal processing and optics.

 \begin{table}[htbp]\label{table:1}
  \centering
  \caption{Some special cases of the multi-dimensional SAFT}
    \begin{tabular}{|c|c|}
    \hline
        Parameters Matrix    &  Corresponding transformation \\
    \hline
    SAFT Matrix $M=\left(\begin{array}{ccc}A&B&P\\C&D&Q\\\end{array}\right)$   & SAFT \\
    \hline
    $A=D=P=Q=\mathbf{0},\ B=-C=I_{N}$ &  FT\\
    $P=Q=\mathbf{0}$ &  LCT\\
    $ A=diag(a_{11},\cdots,a_{NN}),\ B=diag(b_{11},\cdots,b_{NN}),$ & separable LCT  \\
    $C=diag(c_{11},\cdots,c_{NN}),\ D=diag(d_{11},\cdots,d_{NN}),\ P=Q=0$   &\\
    $A=D=diag(\cos\theta_{1},\cdots,\cos\theta_{N}),$  & separable FRFI \\
    $B=-C=diag(\sin\theta_{1},\cdots,\sin\theta_{N}),\ P=Q=0$  & \\
    $A=D=I_{N},\ B=(b_{kl}),\  C=P=Q=\mathbf{0}$  & nonseparable FRT \\
    $A=D=I_{N},\ B=diag(b_{11},\cdots,b_{NN}),\ C=P=Q=\mathbf{0}$ & separable FRT \\
    $A=D=diag(\cosh\phi_{1},\cdots,\cosh\phi_{N}),$  & separable LOT \\
    $B=C=diag(\sinh\phi_{1},\cdots,\sinh\phi_{N}),\ P=Q=0$& \\
\hline
    \end{tabular}%
  \label{tab:addlabel}%
\end{table}%

As we known,  sampling and reconstruction theory
 bridges the modern digital world and the analogue world
of continuous functions. In many applications, taking samples on an appropriate sampling set $X$ is not
practical or even possible. While the dynamical sampling refers to not only the signal $f$
that is sampled but also its various states at different times $(\{t_{0},t_{1},\cdots,t_{n}\})$.

The dynamical sampling problem can be stated as follow: Assume that the initial
state of a physical or a signal process is a function $f$. The operators $A_{t}$
is an evolution operator indexed by $t>0$  and $A_{t}f$ is the
state of $f$ at time $t$. Now we want to recover $f$ from the space samples $\{f(X),f_{t_{1}}(X),\cdots,f_{t_{n}}(X)\}$ of $f$  on $X\in \Omega$, where $\Omega$ is the
domain of $f$, and its various time states $f_{t}(X):=(A_{t}f)(X)$ at times  $t_{0},t_{1},\cdots,t_{n}$. Compared with classical
sampling techniques,  dynamic sampling has potential applications to wireless sensor networks in the health, environment/weather monitoring, ecology,  precision agriculture industries [19-24] and many other fields [21].

At present, Aldroubi, Davis and other researchers have been studied the dynamic sampling for shift-invariant space \cite{aj,ra}.
  With applications in multi-dimensional signal
systems, Zhang et al \cite{zhang} considered the multivariate dynamical sampling in shift-invariant spaces associated with the d-dimensional FT.

To the best of our knowledge, there have been no studies on the multivariate dynamical sampling associated with the multi-dimensional SAFT.
Motivated by \cite{zhang},  we extended the dynamical sampling results associated with the FT or LCT to the multi-dimensional SAFT transform   and shift-invariant space  in the SAFT domains.

Our paper is organized as follows.  Section 2 introduces basic concepts and definitions of the  multi-dimensional SAFT, and its inverse transform.  In Section 3, we consider  convolution operations associated with the SAFT, one for functions and one for
a sequence of numbers and a function. In Section 4,  we derive a sufficient and necessary
condition for sequences in shift invariant spaces under which they can be recovered from
their dynamical sampling measurements in a stable way associated with the multi-dimensional SAFT. In Section
5, we present an example to elucidate our main result, and conclude this
paper in Section 6.

\section{Basic Notations and Definitions}

\ \ \ \ \ Let us give some basic notations and definitions used in this paper.
\begin{definition}
 For any pair functions $f, g\in L^{2}(\mathbb{R}^{2})$, the  convolution operator $*_{c}$   associated with the nD-SAFT  is defined as
\begin{equation}\label{fg}
(f*_{c}g)(t)=\frac{\overline{\lambda}_\mathcal{M}(t)}{\sqrt{|\mathrm{det}(B)|}}(\mathop{f}\limits^{\rightarrow}*\mathop{g}\limits^{\rightarrow})(t),
\end{equation}
where $*$ denotes the traditional convolution operator, $\overline{z}$ means the conjugate of $z$, i.e.,
\begin{equation}\label{tfg}
(f*g)(t)=\int_{\mathbb{R}^{n}}f(t-x)g(x)dx.
\end{equation}
\end{definition}

\begin{definition}
Let $s=s(k)$ be a sequence in $l^{2}(\mathbb{Z}^{n})$, i.e.,$\sum_{k\in \mathbb{Z}^{n}}|s(k)|^{2}\leq +\infty$.
The nD-DT-SAFT of $s$ is defined by
\begin{equation}\label{ss}
(S_{\mathcal{M}}s)(w)
=\frac{1}{\sqrt{|\mathrm{det}(B)|}}\sum_{k\in \mathbb{Z}^{n}}s(k)\lambda_{\mathcal{M}}(k)\eta_{\mathcal{M}}(w)e^{2i\pi[(B^{-1}P)^{T}k-k^{T}B^{-1}w]}.
\end{equation}
\end{definition}

\begin{definition}
The  convolution operator  of a sequence $s\in l^{2}(\mathbb{Z}^{n})$ and a function $\phi\in L^{2}(\mathbb{R}^{n})$ is defined as
\begin{equation}\label{sp}
h(t)=(s*_{sd}\phi)(t)=\frac{\overline{\lambda}_{\mathcal{M}}(t)}{\sqrt{|\mathrm{det}(B)|}}\sum_{k\in \mathbb{Z}^{n}}\lambda_{\mathcal{M}}(k)s(k)\lambda_{\mathcal{M}}(t-k)\phi(t-k).
\end{equation}
\end{definition}

\begin{definition}
A function $\phi(t)$ is said to be periodic with periodicity matrix $M\in\mathbb{R}^{n\times n}$,if
\begin{equation}\label{pm}
\phi(t+Mn)=\phi(t),\ \ for\ \ all \ \ t\in \mathbb{R}^{n}\ \ and \ \ n\in \mathbb{Z}^{n},
\end{equation}
where $M$ is non-singular.
\end{definition}

\begin{definition}\cite{NGY}
The n-dimensional Special Affine Fourier Transform (nD-SAFT) with parameter
 \begin{equation*}
 \mathcal{M}=
\begin{pmatrix}
\begin{tabular}{cc|c}
   A&B&P\\
  C&D&Q\\
\end{tabular}
\end{pmatrix}
 \end{equation*}
 of the signal $f(t)\in L^{2}(\mathbb{R}^{n})$, is defined as
 \begin{equation}\label{s1}
 \begin{aligned}
(S_{\mathcal{M}}f)(w)=\frac{1}{\sqrt{|\mathrm{det}(B)}|}\int_{\mathbb{R}^{n}}f(t)e^{i\pi(t^{T}B^{-1}At+w^{T}DB^{-1}w-2t^{T}B^{-1}w)}\\
\times e^{2i\pi(B^{-1}P)^{T}t+2i\pi(Q^{T}-P^{T}DB^{-1})w}dt,
\end{aligned}
\end{equation}
where $t, w, P, Q\in\mathbb{R}^{n}$ are four real column vectors, and $A,B,C,D$ are $n\times n$ real matrices with  $B$ being non-singular.
The matrix $\mathcal{M}$ is real and  the following equations hold:
\begin{equation*}
AB^{T}=BA^{T},  CD^{T}=DC^{T},  AD^{T}-BC^{T}=I,
\end{equation*}
or
\begin{equation*}
A^{T}C=C^{T}A,  B^{T}D=D^{T}B,  A^{T}D-C^{T}B=I.
\end{equation*}
\end{definition}

nD-SAFT includes several known transforms as special cases. For example, for $P=Q=0$, we get the n-dimensional LCT(\cite{zhang}), which is given by
 \begin{equation}\label{lf}
L\{f(t)\}(w)=\frac{1}{\sqrt{|\mathrm{det}(B)|}}\int_{\mathbb{R}^{n}}f(t)e^{i\pi(t^{T}B^{-1}At+w^{T}DB^{-1}w-2t^{T}B^{-1}w)}dt.
\end{equation}

Since $u$ in \cite{zhang} is a row vector and $t$ in our paper is a column vector, i.e. $u^{T}=w,x=t^{T}$, these two equations are completely identical.

The corresponding inverse formula of nD-LCT $f(x)=L_{M^{-1}}(F)$ is given by
 where $M^{-1}=\begin{pmatrix}
   D^{T}&-B^{T}\\
    -C^{T}& A^{T}
\end{pmatrix}$.
When $A=D=P=Q=0, B =-I, C = I$, we obtain
the Fourier transform. Our goal of this article is to extend these results
to the nD-SAFT domain. In the following
 we consider the inversion formula for the nD-SAFT domain.

\begin{lemma}
The inversion formula for the nD-SAFT is  shown to be
\begin{equation*}
\begin{aligned}
(S_{\mathcal{M}^{-1}}F)(t)
=&\frac{1}{\sqrt{|\mathrm{det}(B)}|}\int_{\mathbb{R}^{n}}F(w)e^{-i\pi[w^{T}(DB^{-1})^{T}w+t^{T}(B^{-1}A)^{T}t-2w^{T}(B^{-1})^{T}t]}\\
&\times e^{-2i\pi(Q^{T}-P^{T}DB^{-1})w-2i\pi(B^{-1}P)^{T}t}dw,
\end{aligned}
\end{equation*}
which may be considered as the nD-SAFT evaluated using the notation $\mathcal{M}^{-1}$ where
\begin{equation*}
 \mathcal{M}^{-1}:=
 \begin{pmatrix}
 D^{T}& -B^{T}& B^{T}Q-D^{T}P\\
  -C^{T}& A^{T}& C^{T}P-A^{T}Q
\end{pmatrix}.
 \end{equation*}
\end{lemma}
\begin{proof}
From \cite{zhang}, we know
\begin{equation}\label{lft}
L^{-1}\{F(w)\}(t)=\frac{1}{\sqrt{|\mathrm{det}(B)}|}\int_{\mathbb{R}^{n}}F(w)e^{-i\pi[w^{T}(DB^{-1})^{T}w+t^{T}(B^{-1}A)^{T}t-2w^{T}(B^{-1})^{T}t]}dw.
\end{equation}
From (\ref{s1}) and (\ref{lf}), we can easily derive the
connection between the nD-SAFT and the nD-LCT as follows:
\begin{equation*}
L\{e^{2i\pi(B^{-1}P)^{T}t}f(t)\}(w)=e^{-2i\pi(Q^{T}-P^{T}DB^{-1})w}(S_\mathcal{M}f)(w).
\end{equation*}
Then, do inverse transformation of nD-LCT on both sides of the above equation
\begin{equation*}
f(t)=e^{-2i\pi(B^{-1}P)^{T}t}L^{-1}[e^{-2i\pi(Q^{T}-P^{T}DB^{-1})w}(S_\mathcal{M}f)(w)].
\end{equation*}
Therefore, using (\ref{lft}), we have
\begin{equation*}
 \begin{aligned}
f(t)
=&\frac{1}{\sqrt{|\mathrm{det}(B)}|}\int_{\mathbb{R}^{n}}e^{-2i\pi(B^{-1}P)^{T}t}e^{-2i\pi(Q^{T}-P^{T}DB^{-1})w}(S_\mathcal{M}f)(w)\\
&\times e^{-i\pi[w^{T}(DB^{-1})^{T}w+t^{T}(B^{-1}A)^{T}t-2w^{T}(B^{-1})^{T}t]}dw\\
=&S_{\mathcal{M}^{-1}}(S_{\mathcal{M}}f).
\end{aligned}
\end{equation*}
This completes the proof.
\end{proof}

For convenience, we introduce notations
\begin{equation*}
\eta_{\mathcal{M}}(w)=e^{i\pi w^{T}DB^{-1}w+2i\pi(Q^{T}-P^{T}DB^{-1})w},
\end{equation*}
and
 $\lambda_{\mathcal{M}}(t)=e^{i\pi t^{T}B^{-1}At}$ is the chirp-modulation
function, so
(\ref{s1}) can be written as
\begin{equation}\label{s2}
(S_\mathcal{M}f)(w)=\frac{1}{\sqrt{|\mathrm{det}(B)}|}\int_{\mathbb{R}^{n}}f(t)\lambda_\mathcal{M}(t)\eta_\mathcal{M}(w)e^{2i\pi[(B^{-1}P)^{T}t-t^{T}B^{-1}w]}dt
\end{equation}

\section{Convolution theorems}
\hskip \parindent
In this section, we study the convolution theorem in the nD-SAFT domain. First, we introduce the notations
\begin{equation*}
\mathop{f}\limits^{\rightarrow}(t):=\lambda_{\mathcal{M}}(t)f(t)=e^{i\pi t^{T}B^{-1}At}f(t)
\end{equation*}

 We will show a convolution theorem for the nD-SAFT as follows.

\begin{prop}\cite{NGY}
For two functions $f, g\in L^{2}(\mathbb{R}^{n})$, let $h(t)=(f*_{c}g)(t)$. Then,
\begin{equation}\label{sfg}
(\mathcal{S}_{\mathcal{M}}h)(w)=\overline{\eta}_{\mathcal{M}}(w)(S_{\mathcal{M}}f)(w)(S_{\mathcal{M}}g)(w).\end{equation}
\end{prop}

We next give the  convolution of a sequence in $l^{2}(\mathbb{Z}^{n})$ and a function in $L^{2}(\mathbb{R}^{n})$, and periodic functions, respectively.

\begin{theorem}
Let $s$ and $\phi$ be as above and $h(t)=(s*_{sd}\phi)(t)$.
Then
\begin{equation}\label{ssp}
(S_{\mathcal{M}}h)(w)=\overline{\eta}_{\mathcal{M}}(w)(S_{\mathcal{M}}s)(w)(S_{\mathcal{M}}\phi)(w),
\end{equation}
Moreover, $|(S_{\mathcal{M}}s)(w)|$  is periodic with periodicity matrix $B$.
\end{theorem}

\begin{proof}
According to (4), (\ref{ss}) and (\ref{sp}), we have
\begin{align*}
&(S_{\mathcal{M}}h)(w)\\
\overset{(4)}=&\frac{1}{\sqrt{|\mathrm{det}(B)|}}\int_{\mathbb{R}^{n}}(s*_{sd}\phi)(t)\lambda_{\mathcal{M}}(t)\eta_{\mathcal{M}}(w)e^{2i\pi[(B^{-1}P)^{T}t-t^{T}B^{-1}w]}dt\\
\overset{(\ref{sp})}=&\frac{1}{\sqrt{|\mathrm{det}(B)|}}\sum_{k\in \mathbb{Z}^{n}}\lambda_{\mathcal{M}}(k)s(k)\int_{\mathbb{R}^{n}}\frac{\overline{\lambda}_{\mathcal{M}}(t)}{\sqrt{|\mathrm{det}(B)|}}\lambda_{\mathcal{M}}(t-k)\phi(t-k)\\
&\times \lambda_{\mathcal{M}}(t)\eta_{\mathcal{M}}(w)e^{2i\pi[(B^{-1}P)^{T}t-t^{T}B^{-1}w]}dt\\
  =&\frac{1}{|\mathrm{det}(B)|}\sum_{k\in \mathbb{Z}^{n}}s(k)\lambda_{\mathcal{M}}(k)\eta_{\mathcal{M}}(w)\int_{\mathbb{R}^{n}}\lambda_{\mathcal{M}}(t-k)\phi(t-k)e^{2i\pi[(B^{-1}P)^{T}t-t^{T}B^{-1}w]}dt\\
  \overset{t-k=t_{1}}{=}&\frac{1}{|\mathrm{det}(B)|}\sum_{k\in \mathbb{Z}^{n}}s(k)\lambda_{\mathcal{M}}(k)\eta_{\mathcal{M}}(w)\int_{\mathbb{R}^{n}}\lambda_{\mathcal{M}}(t_{1})\phi(t_{1})
 e^{2i\pi[
 (B^{-1}P)^{T}(t_{1}+k)-(t_{1}+k)^{T}B^{-1}w]}dt_{1}\\
 =&\frac{1}{\sqrt{|\mathrm{det}(B)|}}\sum_{k\in \mathbb{Z}^{n}}s(k)\lambda_{\mathcal{M}}(k)\eta_{\mathcal{M}}(w)e^{2i\pi[(B^{-1}P)^{T}k-k^{T}B^{-1}w]}\\
 &\times\frac{1}{\sqrt{|\mathrm{det}(B)|}}\int_{\mathbb{R}^{n}}\phi(t_{1})\lambda_{\mathcal{M}}(t_{1})
 e^{2i\pi [(B^{-1}P)^{T}t_{1}-t_{1}^{T}B^{-1}w]
}dt_{1}\\
\overset{(\ref{ss})}=&\overline{\eta}_{\mathcal{M}}(w)(S_{\mathcal{M}}s)(w)(S_{\mathcal{M}}\phi)(w).
\end{align*}
Furthermore, for all $w\in \mathbb{Z}^{n}$ and any $l\in \mathbb{Z}^{n}$, we have
\begin{align*}
&(S_{\mathcal{M}}s)(w+Bl)\\
=&\frac{1}{\sqrt{|\mathrm{det}(B)|}}\sum_{k\in \mathbb{Z}^{n}}s(k)\lambda_{\mathcal{M}}(k)\eta_{\mathcal{M}}(w+Bl)e^{2i\pi[(B^{-1}P)^{T}k-2k^{T}B^{-1}(w+Bl)]}\\
=&\frac{1}{\sqrt{|\mathrm{det}(B)|}}\sum_{k\in \mathbb{Z}^{n}}s(k)\lambda_{\mathcal{M}}(k)e^{i\pi [(w+Bl)^{T}DB^{-1}(w+Bl)+2(Q^{T}-P^{T}DB^{-1})(w+Bl)]}\\
&\times e^{2i\pi[(B^{-1}P)^{T}k-2k^{T}B^{-1}(w+Bl)]}\\
=&\frac{1}{\sqrt{|\mathrm{det}(B)|}}\sum_{k\in \mathbb{Z}^{n}}s(k)\lambda_{\mathcal{M}}(k)\eta_{\mathcal{M}}(w)e^{2i\pi[(B^{-1}P)^{T}k-2k^{T}B^{-1}(w)]}\\
&\times e^{i\pi [(Bl)^{T}DB^{-1}w+ (Bl)^{T}Dl+ w^{T}Dl-2 k^{T}l+2(Q^{T}-P^{T}DB^{-1})Bl]}\\
=&(S_{\mathcal{M}}s)(w)e^{i\pi [(Bl)^{T}DB^{-1}w+ (Bl)^{T}Dl+ w^{T}Dl+2(Q^{T}-P^{T}DB^{-1})Bl]},
\end{align*}
where we use $e^{-2i\pi k^{T}l}=1,(k\in \mathbb{Z}^{n})$ in the last step. Hence,
\begin{equation*}
|(S_{\mathcal{M}}s)(w+Bl)|=|(S_{\mathcal{M}}s)(w)|.
\end{equation*}
This completes the proof.
\end{proof}

Similarly, we introduce the definition and property of the canonical convolution of two sequences in $ l^{2}(\mathbb{Z}^{n})$ associated with the nD-NT-SAFT as
follows.
\begin{definition}
 Let $s=s(k)\in l^{2}(\mathbb{Z}^{n})$ and $c=c(k)\in l^{2}(\mathbb{Z}^{n})$. The canonical convolution operator $*_{d}$ of two sequences $s$ and $c$ is defined by
 \begin{equation}\label{lsc}
 h(l)=(s*_{d}c)(l)=\frac{\overline{\lambda}_{\mathcal{M}}(l)}{\sqrt{|\mathrm{det}(B)|}}\sum_{k\in \mathbb{Z}^{n}}\lambda_{\mathcal{M}}(k)s(k)\lambda_{\mathcal{M}}(l-k)c(l-k).
 \end{equation}
 \end{definition}

\begin{theorem}\label{smhw}
Let $s$ and $c$ be as above and $h(t)=(s*_{d}c)(l)$.
Then
\begin{equation}\label{ssc}
(S_{\mathcal{M}}h)(w)=\overline{\eta}_{\mathcal{M}}(w)(S_{\mathcal{M}}s)(w)(S_{\mathcal{M}}c)(w).
\end{equation}
\end{theorem}

\begin{proof}
By (\ref{sp}) and (\ref{lsc}), we have
\begin{align*}
&(S_{\mathcal{M}}h)(w)\\
=&\frac{1}{\sqrt{|\mathrm{det}(B)|}}\sum_{l\in\mathbb{Z}^{n}}(s*_{d}c)(l)\lambda_{\mathcal{M}}(l)\eta_{\mathcal{M}}(w)e^{2i\pi[(B^{-1}P)^{T}l-l^{T}B^{-1}w]}\\
\overset{(\ref{lsc})}=&\frac{1}{\sqrt{|\mathrm{det}(B)|}}\sum_{l\in \mathbb{Z}^{n}}\lambda_{\mathcal{M}}(l)\eta_{\mathcal{M}}(w)e^{2i\pi[(B^{-1}P)^{T}l-l^{T}B^{-1}w]}\\
&\times\sum_{k\in \mathbb{Z}^{n}}\frac{\overline{\lambda}_{\mathcal{M}}(l)}{\sqrt{|\mathrm{det}(B)|}}\lambda_{\mathcal{M}}(k)s(k)\lambda_{\mathcal{M}}(l-k)c(l-k)\\
  =&\frac{1}{|\mathrm{det}(B)|}\sum_{k\in \mathbb{Z}^{n}}s(k)\lambda_{\mathcal{M}}(k)\eta_{\mathcal{M}}(w)\sum_{l\in \mathbb{Z}^{n}}c(l-k)\lambda_{\mathcal{M}}(l-k)e^{2i\pi[(B^{-1}P)^{T}l-l^{T}B^{-1}w]}\\
 =&\frac{1}{\sqrt{|\mathrm{det}(B)|}}\sum_{k\in \mathbb{Z}^{n}}s(k)\lambda_{\mathcal{M}}(k)\eta_{\mathcal{M}}(w)e^{2i\pi[(B^{-1}P)^{T}k-k^{T}B^{-1}w]}\\
 &\times\frac{\overline{\eta}_{\mathcal{M}}(w)}{\sqrt{|\mathrm{det}(B)|}}\sum_{l\in \mathbb{Z}^{n}}c(l-k)\lambda_{\mathcal{M}}(l-k)\eta_{\mathcal{M}}(w)e^{2i\pi[(B^{-1}P)^{T}(l-k)-(l-k)^{T}B^{-1}w]}\\
  \overset{l-k=l_{1}}{=}&\frac{1}{\sqrt{|\mathrm{det}(B)|}}\sum_{k\in \mathbb{Z}^{2}}s(k)\lambda_{\mathcal{M}}(k)\eta_{\mathcal{M}}(w)e^{2i\pi[(B^{-1}P)^{T}k-k^{T}B^{-1}w]}\\
 &\times\frac{\overline{\eta}_{\mathcal{M}}(w)}{\sqrt{|\mathrm{det}(B)|}}\sum_{l\in \mathbb{Z}^{n}}c(l_{1})\lambda_{\mathcal{M}}(l_{1})\eta_{\mathcal{M}}(w)e^{2i\pi[(B^{-1}P)^{T}l_{1}-l_{1}^{T}B^{-1}w]}\\
 =&\overline{\eta}_{\mathcal{M}}(w)(S_{\mathcal{M}}s)(w)(S_{\mathcal{M}}c)(w).
\end{align*}
This completes the proof.

\end{proof}

 On the basis of the above lemmas, we propose to explore the multivariate
dynamical sampling in shift-invariant spaces associated with the nD-NS-SAFT in
what follows.

\section{ Multivariate dynamical sampling in shift-invariant spaces }
\hskip \parindent
In this section, we naturally study a shift-invariant spaces associated with the 2D-SAFT.  we give a necessary and sufficient condition for a function $\phi(t)\in L^{2}(\mathbb{R}^{n})$ to be a generator for
a shift-invariant space in terms of its SAFT.\\
Let $M$ be a $n\times n$ non-singular real matrix (not necessarily integer matrix)
with $m=|\mathrm{det}(M)|$, and $MZ^{n} \overset{\triangle}{=} \{Mn | n \in\mathbb{Z}^{n}\}$ be a lattice generated by
$M$ [7, 8]. Let us define $\mathbb{T}^{n}\overset{\triangle}{=}\{[x_{1}, x_{2},\cdots,x_{n}]^{T}
| x_{i}\in[0, 1) \} \subset \mathbb{R}^{n}$.
The fundamental parallelepiped of $MZ^{n}$
is defined as the region
\begin{equation*}
\mathcal{G}(M)\overset{\triangle}{=}\{Mx|x\in\mathbb{T}^{n}\}.
\end{equation*}
One can see that $\mathcal{G}(M)$ and its shifted copies (called the other lattice cells)
constitute the whole real vector space $\mathbb{R}^{n}$
, i.e.,
\begin{equation}\label{bc}
\bigcup_{n\in\mathbb{Z}^{n}}\{M(x+n)|x\in\mathbb{T}^{n}\}=\mathbb{R}^{n}.
\end{equation}

When $M$ is further an integer matrix, we define
\begin{equation*}
\mathcal{N}(M)\overset{\triangle}{=}\{k|k=Mx, x\in\mathbb{T}^{n}, \ \ and\ \ k\in \mathbb{Z}^{n}\}.
\end{equation*}
 The number of elements in $\mathcal{N}(M)$ equals $m$. Without loss of generality, let $\gamma_{0}=[0,0]^{T}, \gamma_{1},\dots$,
 $\gamma_{m-1}$ be the $m$ distinct elements in $\mathcal{N}(M)$. It is clear that $\gamma_{j}+M\mathbb{Z}^{n}(j=0,1,\dots,m-1)$ constitute the whole integer vector space $\mathbb{Z}^{n}$, i.e., $U_{j=0}^{m-1}\{\gamma_{j}+M\mathbb{Z}^{n}\}=\mathbb{Z}^{n}$.  Analogously, let $M^{T}\mathbb{Z}^{n}\overset{\triangle}{=}\{M^{T}n|n\in\mathbb{Z}^{n}\}$ be a lattice generated by $M^{T}$, and $\eta_{0}=[0,0]^{T}, \eta_{1},\dots,\eta_{m-1}$
be the $m$ distinct elements in $\mathcal{N}(M^{T})$.
Obviously, $U_{j=0}^{m-1}\{\eta_{j}+M^{T}\mathbb{Z}^{n}\}=\mathbb{Z}^{n}$.

\subsection{ Shift-invariant spaces }
\hskip \parindent
\begin{theorem}
Let a sequence ${s(n)}\in l^{2}(\mathbb{R}^{n})$,and a function $\phi(t)\in L^{2}(\mathbb{R}^{n})$, Assume that the chirp-modulated subspace of $L^{2}(\mathbb{R}^{n})$
 is given by
 \begin{equation*}
 V(\phi)=\overline{\{f\in L^{2}(\mathbb{R}^{n}):f(t)=(s*_{sd}\phi)(t)\}}.
 \end{equation*}
 Then $\left\{e^{-2i\pi (t-k)^{T}B^{-1}(t-k)}\phi(t-k)\right\}$ is a Riesz basis for $V(\phi)$, if and only if
there exist two constants $\eta_{1}, \eta_{2} > 0$ such that
\begin{equation}\label{eta}
\eta_{1}\leq\sum_{k=-\infty}^{+\infty}|(S_{\mathcal{M}}\phi)(w+Bk)|\leq\eta_{2}
\end{equation}
for all $w\in\mathcal{G}(B)=\{Bx|x\in \mathbb{T}^{n}\}$.
\end{theorem}
\begin{proof}
For $f(t)=(s*_{sd}\phi)(t)$, by Theorem 3.6, we have
\begin{equation*}
(S_{\mathcal{M}}f)(w)=\overline{\eta}_{\mathcal{M}}(w)(S_{\mathcal{M}}s)(w)(S_{\mathcal{M}}\phi)(w).
\end{equation*}
By using $|\overline{\eta}_{\mathcal{M}}(w)|=1$,
\begin{equation*}
|(S_{\mathcal{M}}f)(w)|^{2}=|(S_{\mathcal{M}}s)(w)|^{2}|(S_{\mathcal{M}}\phi)(w)|^{2}.
\end{equation*}
Since $|(S_{\mathcal{M}}s)(w)|$ is periodic with periodicity matrix $B$, we have, from (\ref{bc}),
\begin{equation}\label{SG}
\begin{aligned}
&\Vert(S_{\mathcal{M}}f)(w)\Vert^{2}_{L^{2}(\mathbb{R}^{n})}\\
=&\int_{\mathbb{R}^{n}}|(S_{\mathcal{M}}s)(w)|^{2}|(S_{\mathcal{M}}\phi)(w)|^{2}dw\\
=&\sum_{k\in\mathbb{Z}^{n}}\int_{\{Bk+B\mathbb{T}^{n}\}}|(S_{\mathcal{M}}s)(w)|^{2}|(S_{\mathcal{M}}\phi)(w)|^{2}dw\\
\overset{w=w_{1}+Bk}{=}&\sum_{k\in\mathbb{Z}^{n}}\int_{\{B\mathbb{T}^{n}\}}|(S_{\mathcal{M}}s)(w_{1}+Bk)|^{2}|(S_{\mathcal{M}}\phi)(w_{1}+Bk)|^{2}dw_{1}\\
=&\int_{\{B\mathbb{T}^{2}\}}|(S_{\mathcal{M}}s)w_{1}|^{2}\sum_{k\in\mathbb{Z}^{n}}|(S_{\mathcal{M}}\phi)(w_{1}+Bk)|^{2}dw_{1}\\
:=&\int_{\{B\mathbb{T}^{n}\}}|(S_{\mathcal{M}}s)w|^{2}G_{M,\phi}(w)dw,
\end{aligned}
\end{equation}
where $G_{M,\phi}(w)=\sum_{k\in\mathbb{Z}^{n}}|(S_{\mathcal{M}}\phi)(w+Bk)|^{2}$ is the Grammian of $\phi$ associated
with the nD-SAFT. Notice that

\begin{align*}
&\int_{\{Bt|t\in\mathbb{T}^{n}\}}|(S_{\mathcal{M}}s)(w)|^{2}dw\\
=&\int_{\{B\mathbb{T}^{n}\}}\frac{1}{\sqrt{|\mathrm{det}(B)|}}\sum_{n\in \mathbb{Z}^{n}}s(n)\lambda_{\mathcal{M}}(n)\eta_{\mathcal{M}}(w)e^{i\pi[(B^{-1}P)^{T}n-2n^{T}B^{-1}w]}\\
&\times\frac{1}{\sqrt{|\mathrm{det}(B)|}}\sum_{k\in \mathbb{Z}^{n}}\overline{s(k)}\overline{\lambda}_{\mathcal{M}}(k)\overline{\eta}_{\mathcal{M}}(w)e^{i\pi[(B^{-1}P)^{T}k-2k^{T}B^{-1}w]}\\
=&\frac{1}{|\mathrm{det}(B)|}\sum_{n\in \mathbb{Z}^{n}}\sum_{k\in \mathbb{Z}^{n}}s(n)\overline{s(k)}\lambda_{\mathcal{M}}(n)\overline{\lambda}_{\mathcal{M}}(k)e^{i\pi(B^{-1}P)^{T}(n-k)}\\
&\times\int_{\{B\mathbb{T}^{n}\}}e^{-2i\pi(n-k)^{T}B^{-1}w}dw,\\
\end{align*}

and since
\begin{equation*}
\int_{\{B\mathbb{T}^{n}\}}e^{-2i\pi(n-k)^{T}B^{-1}w}dw=
|\mathrm{det}(B)|\delta_{n,k},
\end{equation*}
it follows that
\begin{equation}\label{sLl}
\begin{aligned}
\int_{\{Bt|t\in\mathbb{T}^{n}\}}|(S_{\mathcal{M}}s)(w)|^{2}dw=\sum_{n\in \mathbb{Z}^{n}}|s(n)|^{2}=||s||^{2}_{l^{2}(\mathbb{R}^{n})}.
\end{aligned}
\end{equation}
From (\ref{sLl}), we know that
\begin{equation*}
0<\eta_{1}\leq G_{M,\phi}(w)\leq\eta_{2}<+\infty
\end{equation*}
is equivalent to
\begin{equation*}
\eta_{1}\left\Vert(S_{\mathcal{M}}s)(w)\right\Vert^{2}=\eta_{1}||s(k)||^{2}_{l^{2}}\leq \left\Vert(S_{\mathcal{M}}f)(w)\right\Vert^{2}_{L^{2}(\mathbb{R}^{n})}\leq \eta_{2}||s(k)||^{2}_{l^{2}}\leq \eta_{2}\left\Vert(S_{\mathcal{M}}s)(w)\right\Vert^{2}.
\end{equation*}
This completes the proof.
\end{proof}

The local behavior and global decay of $\phi$ can be described in terms of the Wiener amalgam spaces as follows. A measurable function $f$ belongs to the Wiener amalgam space $W(L^{p}(\mathbb{R}^{n})),1\leq p<+\infty$, if it satisfies
\begin{equation}\label{wnas}
||f||^{p}_{W(L^{p}(\mathbb{R}^{n}))}:=
\sum_{k\in{\mathbb{Z}^{n}}}ess sup\{|f(x+k)|^{p};x\in \mathbb{T}^{n}\}<+\infty.
\end{equation}
If $p=\infty$,  a measurable function $f$ belongs to $W(L^{\infty}(\mathbb{R}^{n}))=L^{\infty}(\mathbb{R}^{n})$, if it satisfies
\begin{equation}
||f||^{p}_{W(L^{\infty}(\mathbb{R}^{n}))}:=
\sum_{k\in{\mathbb{Z}^{n}}}ess sup\{|f(x+k)|;x\in \mathbb{T}^{n}\}<+\infty.
\end{equation}

Considering that ideal sampling makes sense only for continuous functions, we therefore focus on the amalgam space $W_{0}(L^{p}(\mathbb{R}^{n})):=W(L^{p}(\mathbb{R}^{n}))\cap C(\mathbb{R}^{n})$, where $C(\mathbb{R}^{n})$  denotes the space of continuous functions on $\mathbb{R}^{n}$.

The multivariate dynamical sampling problem in shift-invariant spaces is
to recover a function $f\in V(\phi)$  from its dynamical sampling measurements
associated with the 2D-NS-SAFT, i.e.,
\begin{equation*}
\{(a^{j}*_{c}f)(M^{T}k):j=1, 2, \cdots, m-1, k\in\mathbb{Z}^{n}\},
\end{equation*}
where $a^{j}=a*_{c}\cdots_{c}a$ and $a\in W(L^{1}(\mathbb{R}^{n}))$.

First, we present an important lemma, which will be used later.
\begin{lemma}
Let a sequence $s\in l^{2}(\mathbb{Z}^{n})$ and two functions $f,g\in L^{2}(\mathbb{R}^{n})$. Then, we have
\begin{equation}\label{fsg}
f*_{c}(s*_{sd}g)=s*_{sd}(f*_{c}g).
\end{equation}
\end{lemma}

\begin{proof}
Letting $u_{1}(t)=(s*_{sd}g)(t)$ and $u_{2}(t)=(f*_{c}g)(t)$,
we have
\begin{align*}
&(f*_{c}u_{1})(t)=\frac{\overline{\lambda}_{\mathcal{M}}(t)}{\sqrt{|\mathrm{det}(B)|}}(\mathop{f}\limits^{\rightarrow}*\mathop{u_{1}}\limits^{\rightarrow})(t)\\
=&\frac{1}{\sqrt{|\mathrm{det}(B)|}}e^{-i\pi t^{T}B^{-1}At}\int_{\mathbb{R}^{n}}f(t-x)e^{i\pi (t-x)^{T}B^{-1}A(t-x)}u_{1}(x)e^{i\pi x^{T}B^{-1}Ax}dx \\
=&\frac{1}{|\mathrm{det}(B)|}e^{-i\pi t^{T}B^{-1}At}\int_{\mathbb{R}^{n}}f(t-x)e^{i\pi (t-x)^{T}B^{-1}A(t-x)}\\
&\times\left(\sum_{k\in \mathbb{Z}^{n}}s(k)e^{i\pi k^{T}B^{-1}Ak}g(x-k)e^{i\pi (x-k)^{T}B^{-1}A(x-k)}\right)dx\\
\overset{y=x-k}{=}&\frac{1}{|\mathrm{det}(B)|}e^{-i\pi t^{T}B^{-1}At}\sum_{k\in \mathbb{Z}^{n}}s(k)e^{i\pi k^{T}B^{-1}Ak}(\int_{\mathbb{R}^{n}}f(t-k-y)\\
&e^{i\pi (t-k-y)^{T}B^{-1}A(t-k-y)}
g(y)e^{i\pi y^{T}B^{-1}Ay}dy)\\
=&\frac{1}{\sqrt{|\mathrm{det}(B)|}}e^{-i\pi t^{T}B^{-1}At}\sum_{k\in \mathbb{Z}^{n}}s(k)e^{i\pi k^{T}B^{-1}Ak}u_{2}(t-k)
e^{i\pi (t-k)^{T}B^{-1}A(t-k)}\\
=&(s*_{sd}u_{2})(t).
\end{align*}

This completes the proof.

\end{proof}
\subsection{ The characterization with the discrete SAFT }
In the following, we propose a sufficient and necessary condition for stably recovering $f$ from its dynamical sampling measurements
$f(M^{T}k),(a^{j}*_{c}f)(M^{T}k), j=1,2,\cdots, m-1, k\in \mathbb{Z}^{2}$ associated with the nD-NS-SAFT.
\begin{theorem}
Let $\phi=\phi_{0}\in W_{0}(L^{1}(\mathbb{R}^{n})), a\in W(L^{1}(\mathbb{R}^{n})$, and $\phi_{j}=a^{j}*_{c}\phi$, then $S_{\mathcal{M}}\phi^{j}_{l}\in C(\mathbb{R}^{n})$ for $l,j=0,1,\cdots, m-1$, where $\phi_{l}^{j}(r)=\phi_{j}(M^{T}r-\eta_{l})e^{i\pi (M^{T}r-\eta_{l})^{T}B^{-1}A(M^{T}r-\eta_{l})}$. Any $f\in V(\phi)$ can be  recovered in a stable
way, i.e., the inverse of $\mathcal{B}(w)$ is bounded, from the dynamical sampling measurements $f(M^{T}k), (a^{j}*_{c}f)(M^{T}k), j=1,2,\cdots, m-1, k\in \mathbb{Z}^{n}$, if and only if $\mathrm{det}(\mathcal{B}(w))\neq 0$ for any $w\in \mathbb{R}^{n}$, where $\mathcal{B}(w)$ is defined by
\begin{equation*}
\mathcal{B}(w)=
\left[{\begin{array}{cccc}
S_{\mathcal{M}}\phi^{0}_{0}(w)&S_{\mathcal{M}}\phi^{0}_{1}(w)&\cdots&S_{\mathcal{M}}\phi^{0}_{m-1}(w)\\
S_{\mathcal{M}}\phi^{1}_{0}(w)&S_{\mathcal{M}}\phi^{1}_{1}(w)&\cdots&S_{\mathcal{M}}\phi^{1}_{m-1}(w)\\
\vdots&\vdots&\vdots&\vdots\\
S_{\mathcal{M}}\phi^{m-1}_{0}(w)&S_{\mathcal{M}}\phi^{m-1}_{1}(w)&\cdots&S_{\mathcal{M}}\phi^{m-1}_{m-1}(w)
\end{array}}\right].
\end{equation*}
\end{theorem}
\begin{proof}
Given a function $f\in V(\phi)$, we have
\begin{equation*}
f=(s*_{sd}\phi)(t).
\end{equation*}
By Lemma 5.2 , we have
\begin{align*}
v_{j}(k)\triangleq &(a^{j}*_{c}f)(M^{T}k)e^{i\pi (M^{T}k)^{T}B^{-1}AM^{T}k-i\pi k^{T}B^{-1}Ak}\\
=&[a^{j}*_{c}(s*_{sd}\phi)](M^{T}k)e^{i\pi (M^{T}k)^{T}B^{-1}AM^{T}k-i\pi k^{T}B^{-1}Ak}\\
\overset{(\ref{fsg})}{=}&[s*_{sd}(a^{j}*_{c}\phi)](M^{T}k)e^{i\pi (M^{T}k)^{T}B^{-1}AM^{T}k-i\pi k^{T}B^{-1}Ak}\\
\overset{(\ref{sp})}{=}&\frac{\overline{\lambda}_{M}(M^{T}k)}{\sqrt{|\mathrm{det}(B)|}}\sum_{n\in \mathbb{Z}^{n}}\lambda_{M}(n)s(n)\lambda_{M}(M^{T}k-n)(a^{j}*_{c}\phi)(M^{T}k-n)
\\
&\times e^{i\pi (M^{T}k)^{T}B^{-1}AM^{T}k-i\pi k^{T}B^{-1}Ak}\\
=& \frac{1}{\sqrt{|\mathrm{det}(B)|}}e^{-i\pi k^{T}B^{-1}Ak}\sum_{n\in \mathbb{Z}^{n}}\lambda_{\mathcal{M}}(n)s(n)\lambda_{\mathcal{M}}(M^{T}k-n)\phi_{j}(M^{T}k-n)\\
=&\frac{1}{\sqrt{|\mathrm{det}(B)|}}e^{-i\pi k^{T}B^{-1}Ak}\sum_{l=0}^{m-1}\sum_{r\in \mathbb{Z}^{n}}\lambda_{\mathcal{M}}(M^{T}r+\eta_{l})s(M^{T}r+\eta_{l})\lambda_{\mathcal{M}}(M^{T}k-M^{T}r-\eta_{l})\phi_{j}(M^{T}k-M^{T}r-\eta_{l})\\
=&\sum_{l=0}^{m-1}(s_{l}*_{sd}\phi_{l}^{j})(k),
\end{align*}
where $s_{l}(r)=s(M^{T}r+\eta_{l})$
Then, by Theorem 3.6, we readily have, for $j=1,2,\cdots, m-1$,
\begin{equation}\label{svj}
(S_{\mathcal{M}}v_{j})(w)=\sum_{l=0}^{m-1}e^{-i\pi w^{T}DB^{-1}w+2i\pi(Q^{T}-P^{T}DB^{-1})w}(S_{\mathcal{M}}s_{l})(w)(S_{\mathcal{M}}\phi_{l}^{j})(w).
\end{equation}
Let
\begin{equation*}
(S_{\mathcal{M}}v)(w)=\left[{\begin{array}{c}
(S_{\mathcal{M}}v_{0})(w)\\
(S_{\mathcal{M}}v_{1})(w)\\
\vdots\\
(S_{\mathcal{M}}v_{m-1})(w)
\end{array}}\right],
\end{equation*}
and
\begin{equation*}
(S_{\mathcal{M}}s)(w)=\left[{\begin{array}{c}
(S_{\mathcal{M}}s_{0})(w)\\
(S_{\mathcal{M}}s_{1})(w)\\
\vdots\\
(S_{\mathcal{M}}s_{m-1})(w)
\end{array}}\right].
\end{equation*}
Hence, from (\ref{svj}), we get
\begin{equation}\label{svb}
(S_{\mathcal{M}}v)(w)=\mathcal{B}(w)(S_{\mathcal{M}}s)(w).
\end{equation}
That is to say, we can solve the equation (\ref{svb}) with respect to $(S_{\mathcal{M}}s)(w)$, if $\mathcal{B}(w)$ is invertible.

\end{proof}

\subsection{ The characterization with the continuous SAFT }
In the following, we characterize the problem of multivariate dynamical
sampling in the shift-invariant space $V(\phi)$ with the continuous
SAFT transform.
Before proving The main theorem, we introduce the Poisson summation formula.

\begin{lemma}\label{g(t)}
Suppose $g(t)\in L^{2}(B\mathbb{T}^{n})$, then
the Poisson summation formula for the nD-SAFT which is given by
\begin{equation*}
\frac{\overline{\eta}_\mathcal{M}(w)}{\sqrt{|\mathrm{det}(B)|}}\sum_{k\in \mathbb{Z}^{n}}g(k)\lambda_{\mathcal{M}}(k)\eta_{\mathcal{M}}(w)e^{2i\pi[(B^{-1}P)^{T}k-k^{T}B^{-1}w]}=\sum_{n\in\mathbb{Z}^{n}}\overline{\eta}_\mathcal{M}(B^{-1}w+n)(S_\mathcal{M}g)(B^{-1}w+n)
\end{equation*}
holds almost everywhere.
\end{lemma}
\begin{proof}
Using the Poisson summation formula for Fourier transformation, we have

\begin{align*}
&\frac{1}{\sqrt{|\mathrm{det}(B)|}}\sum_{k\in \mathbb{Z}^{n}}g(k)\lambda_{\mathcal{M}}(k)\eta_{\mathcal{M}}(w)e^{2i\pi[(B^{-1}P)^{T}k-k^{T}B^{-1}w]}\\
=&\frac{\eta_{\mathcal{M}}(w)}{\sqrt{|\mathrm{det}(B)|}}\sum_{k\in \mathbb{Z}^{n}}g(k)\lambda_{\mathcal{M}}(k)e^{2i\pi[(B^{-1}P)^{T}k-k^{T}B^{-1}w]}\\
=&\frac{\eta_{\mathcal{M}}(w)}{\sqrt{|\mathrm{det}(B)|}}\sum_{k\in \mathbb{Z}^{n}}h(k)e^{-2i\pi k^{T}B^{-1}w}\\
=&\frac{\eta_{\mathcal{M}}(w)}{\sqrt{|\mathrm{det}(B)|}}(\mathcal{F}h)(B^{-1}w)\\
=&\frac{\eta_{\mathcal{M}}(w)}{\sqrt{|\mathrm{det}(B)|}}\sum_{n\in \mathbb{Z}^{n}}(\mathcal{F}h)(B^{-1}w+n)\\
=&\frac{\eta_{\mathcal{M}}(w)}{\sqrt{|\mathrm{det}(B)|}}\sum_{n\in \mathbb{Z}^{n}}\int_{\mathbb{R}^{n}}h(t)e^{-2\pi it(B^{-1}w+n)}dt\\
=&\frac{\eta_{\mathcal{M}}(w)}{\sqrt{|\mathrm{det}(B)|}}\sum_{n\in \mathbb{Z}^{n}}\int_{\mathbb{R}^{n}}g(t)\lambda_{\mathcal{M}}(t)e^{2i\pi(B^{-1}P)^{T}t}e^{-2\pi it(B^{-1}w+n)}dt\\
=&\frac{\eta_{\mathcal{M}}(w)}{\sqrt{|\mathrm{det}(B)|}}\sum_{n\in\mathbb{Z}^{n}}\overline{\eta}_\mathcal{M}(B^{-1}w+n)(S_\mathcal{M}g)(B^{-1}w+n),
\end{align*}

where  $h(k)=g(k)\lambda_{\mathcal{M}}(k)\eta_{\mathcal{M}}(w)e^{2i\pi(B^{-1}P)^{T}k}$, the proof is completed.
\end{proof}

\begin{lemma}\label{c(t)}
Suppose $c(t)\in l^{2}(B\mathbb{T}^{n})$, then
we have
\begin{equation*}
m\overline{\eta}_\mathcal{M}(w)\{\mathcal{S_{M}}[D_{M}c]\}(w)
=\sum_{v=0}^{m-1}\overline{\eta}_\mathcal{M}[BM^{-1}(w+\gamma_{v})](S_\mathcal{M}c)[BM^{-1}(w+\gamma_{v})].
\end{equation*}
\end{lemma}
\begin{proof}
By Proposition 4.1, we have
\begin{equation*}
\begin{aligned}
&\sum_{v=0}^{m-1}\overline{\eta}_{\mathcal{M}}[BM^{-1}(w+\gamma_{v})](S_\mathcal{M}c)[BM^{-1}(w+\gamma_{v})]\\
=&\sum_{v=0}^{m-1}\sum_{n\in\mathbb{Z}^{n}}c(n)\lambda_{\mathcal{M}}(n)\eta_{\mathcal{M}}[BM^{-1}(w+\gamma_{v})]e^{2i\pi[(B^{-1}P)^{T}n-n^{T}B^{-1}BM^{-1}(w+\gamma_{v})]}\\
=&\frac{1}{\sqrt{|\mathrm{det}(B)|}}\sum_{n^{'}\in\mathbb{Z}^{n},n=M^{T}n^{'}+\eta_{0}}c(n)\lambda_\mathcal{M}(n)e^{-2i\pi n^{T}M^{-1}w}e^{2i\pi[(B^{-1}P)^{T}n}\sum_{v=0}^{m-1}\eta_{\mathcal{M}}[BM^{-1}(w+\gamma_{v})]e^{-2i\pi n^{T}M^{-1}\gamma_{v}}\\
&+\frac{1}{\sqrt{|\mathrm{det}(B)|}}\sum_{j=1}^{m-1}\sum_{n^{'}\in\mathbb{Z}^{n},n=M^{T}n^{'}+\eta_{j}}c(n)\lambda_\mathcal{M}(n)e^{-2i\pi n^{T}M^{-1}w}e^{2i\pi[(B^{-1}P)^{T}n}\sum_{v=0}^{m-1}\eta_{\mathcal{M}}[BM^{-1}(w+\gamma_{v})]e^{-2i\pi n^{T}M^{-1}\gamma_{v}}\\
\end{aligned}
\end{equation*}
\begin{equation*}
\begin{aligned}
&\sum_{v=0}^{m-1}\overline{\eta}_{\mathcal{M}}[BM^{-1}(w+\gamma_{v})](S_\mathcal{M}c)[BM^{-1}(w+\gamma_{v})]\\
=&\sum_{v=0}^{m-1}\sum_{n\in\mathbb{Z}^{n}}c(n)\lambda_{\mathcal{M}}(n)e^{2i\pi[(B^{-1}P)^{T}n-n^{T}B^{-1}BM^{-1}(w+\gamma_{v})]}\\
=&\frac{1}{\sqrt{|\mathrm{det}(B)|}}\sum_{n^{'}\in\mathbb{Z}^{n},n=M^{T}n^{'}+\eta_{0}}c(n)\lambda_\mathcal{M}(n)e^{2i\pi[(B^{-1}P)^{T}n- n^{T}M^{-1}w]}\sum_{v=0}^{m-1}e^{-2i\pi n^{T}M^{-1}\gamma_{v}}\\
&+\frac{1}{\sqrt{|\mathrm{det}(B)|}}\sum_{j=1}^{m-1}\sum_{n^{'}\in\mathbb{Z}^{n},n=M^{T}n^{'}+\eta_{j}}c(n)\lambda_\mathcal{M}(n)e^{2i\pi[(B^{-1}P)^{T}n- n^{T}M^{-1}w]}\sum_{v=0}^{m-1}e^{-2i\pi n^{T}M^{-1}\gamma_{v}}\\
=&\frac{m}{\sqrt{|\mathrm{det}(B)|}}\sum_{n^{'}\in\mathbb{Z}^{n}}c(M^{T}n^{'})\lambda_\mathcal{M}(M^{T}n^{'})e^{2i\pi[(B^{-1}P)^{T}M^{T}n^{'}- (M^{T}n^{'})^{T}M^{-1}w]}\\
=&m\overline{\eta}_\mathcal{M}(w)\{\mathcal{S_{M}}[D_{M}c]\}(w).
\end{aligned}
\end{equation*}
This completes the proof.
\end{proof}

Let $\phi_{j}=a^{j}*\phi$, $(S_\mathcal{M})\Phi_{j}=\sum_{n\in\mathbb{Z}^{2}}\overline{\eta}^{2}_\mathcal{M}(w+n)(S_\mathcal{M})\phi_{j}(w+n)$ and
\begin{equation*}
\mathcal{D}(w)=
\left[{\begin{array}{cccc}
(S_\mathcal{M})\Phi_{0}(M^{-1}w)&(S_\mathcal{M})\Phi_{0}(M^{-1}(w+\gamma_{1}))&\cdots&(S_\mathcal{M})\Phi_{0}(M^{-1}(w+\gamma_{m-1}))\\
(S_\mathcal{M})\Phi_{1}(M^{-1}w)&(S_\mathcal{M})\Phi_{1}(M^{-1}(w+\gamma_{1}))&\cdots&(S_\mathcal{M})\Phi_{1}(M^{-1}(w+\gamma_{m-1}))\\
\vdots&\vdots&\vdots&\vdots\\
(S_\mathcal{M})\Phi_{m-1}(M^{-1}w)&(S_\mathcal{M})\Phi_{m-1}(M^{-1}(w+\gamma_{1}))&\cdots&(S_\mathcal{M})\Phi_{m-1}(M^{-1}(w+\gamma_{m-1}))\\
\end{array}}\right].
\end{equation*}
Then we have the following result.
\begin{theorem}
Let $\phi=\phi_{0}\in W_{0}(L^{1}(\mathbb{R}^{n})), a\in W(L^{1}(\mathbb{T}^{n}))$,  then $\Phi_{j}\in C(\mathbb{T}^{n})$ for $l,j=0,1,\cdots, m-1$, For any $f\in V(\phi)$ can be  recovered in a stable
way, i.e., the inverse of $\mathcal{D}(w)$ is bounded, from the dynamical sampling measurements $ (a^{j}*f)(M^{T}k), j=1,2,\cdots, m-1, k\in \mathbb{Z}^{n}$, if and only if $\mathrm{det}(\mathcal{D}(w))\neq 0$ for any $w\in \mathbb{R}^{2}$, here $\mathbb{T}^{n}\simeq [0,1)^{n}$ is a torus.
\end{theorem}

We are now ready to prove the main result.

\begin{proof}
Let $f_{j}=a^{j}*f$, $h_{j}=f_{j}|_{\mathbb{Z}^{2}}$. Using Lemma \ref{g(t)} and
 \ref{c(t)}, for any $j=1,2,\cdots, m-1$, we have
 \begin{align*}
&m\{\mathcal{S_{M}}[D_{M}(h_{j})]\}(w)\\
=&\eta_\mathcal{M}(w)\sum_{v=0}^{m-1}\overline{\eta}_\mathcal{M}[BM^{-1}(w+\gamma_{v})](S_\mathcal{M})(h_{j})[BM^{-1}(w+\gamma_{v})]\\
=&\eta_\mathcal{M}(w)\sum_{v=0}^{m-1}\overline{\eta}_\mathcal{M}[BM^{-1}(w+\gamma_{v})]\eta_\mathcal{M}(w)\sum_{n\in\mathbb{Z}^{n}}\overline{\eta}_\mathcal{M}(B^{-1}(BM^{-1}(w+\gamma_{v}))+n)(S_\mathcal{M}f_{j})(B^{-1}(BM^{-1}(w+\gamma_{v}))+n)\\
=&\eta^{2}_\mathcal{M}(w)\sum_{v=0}^{m-1}\overline{\eta}_\mathcal{M}[BM^{-1}(w+\gamma_{v})]\sum_{n\in\mathbb{Z}^{n}}\overline{\eta}_\mathcal{M}(M^{-1}(w+\gamma_{v})+n)S_\mathcal{M}(a^{j}*(s*_{sd}\phi_{j}))(M^{-1}(w+\gamma_{v})+n)\\
=&\eta^{2}_\mathcal{M}(w)\sum_{v=0}^{m-1}\overline{\eta}_\mathcal{M}[BM^{-1}(w+\gamma_{v})]\sum_{n\in\mathbb{Z}^{n}}\overline{\eta}_\mathcal{M}(M^{-1}(w+\gamma_{v})+n)S_\mathcal{M}(s*_{sd}\phi_{j})(M^{-1}(w+\gamma_{v})+n)\\
=&\eta^{2}_\mathcal{M}(w)\sum_{v=0}^{m-1}\overline{\eta}_\mathcal{M}[BM^{-1}(w+\gamma_{v})]\sum_{n\in\mathbb{Z}^{n}}\overline{\eta}^{2}_\mathcal{M}(M^{-1}(w+\gamma_{v})+n)S_\mathcal{M}s(M^{-1}(w+\gamma_{v})+n)S_\mathcal{M}\phi_{j}(M^{-1}(w+\gamma_{v})+n)\\
=&\eta^{2}_\mathcal{M}(w)\sum_{v=0}^{m-1}\overline{\eta}_\mathcal{M}[BM^{-1}(w+\gamma_{v})]S_\mathcal{M}s(M^{-1}(w+\gamma_{v}))\sum_{n\in\mathbb{Z}^{n}}\overline{\eta}^{2}_\mathcal{M}(M^{-1}(w+\gamma_{v})+n)S_\mathcal{M}\phi_{j}(M^{-1}(w+\gamma_{v})+n)\\
=&\eta^{2}_\mathcal{M}(w)\sum_{v=0}^{m-1}\overline{\eta}_\mathcal{M}[BM^{-1}(w+\gamma_{v})]S_\mathcal{M}s(M^{-1}(w+\gamma_{v}))(S_\mathcal{M})\Phi_{j}(M^{-1}(w+\gamma_{v})),
\end{align*}

 where
\begin{equation}\label{qiuhe}
(S_\mathcal{M})\Phi_{j}(w)=\sum_{n\in\mathbb{Z}^{n}}\overline{\eta}^{2}_\mathcal{M}(w+n)S_\mathcal{M}\phi_{j}(w+n).
\end{equation}
 Define
\begin{equation*}
h(w)=\left[{\begin{array}{c}
m\{\mathcal{S_{M}}[D_{M}(h_{0})]\}(w)\\
m\{\mathcal{S_{M}}[D_{M}(h_{1})]\}(w)\\
\vdots\\
m\{\mathcal{S_{M}}[D_{M}(h_{m-1})]\}(w)
\end{array}}\right].
\end{equation*}
and
\begin{equation}\label{cwg}
\begin{aligned}
\mathcal{C}(w)=\left[{\begin{array}{c}
(S_{\mathcal{M}}s)[BM^{-1}(w)]\\
(S_{\mathcal{M}}s)[BM^{-1}(w+\gamma_{1})]\\
\vdots\\
(S_{\mathcal{M}}s)[BM^{-1}(w+\gamma_{m-1})]
\end{array}}\right].
\end{aligned}
\end{equation}
 Then
$$ h(w)=\mathcal{D}(w)\mathcal{C}(w),$$
 where $\mathcal{C}(w)$ is defined by (\ref{cwg}). We can solve the above Equation  with
respect to $\mathcal{C}(w)$ (which we use to produce $f$) if $\mathcal{D}(w)$ is invertible.

\end{proof}

\section{Example of multivariate dynamical sampling}
\hskip \parindent
In this section, to make the main results obtained above more transparent and more complete, we give a simple example to show that the necessary
and sufficient condition in Theorem 4.6 is feasible.

For simplicity, we take $d=2$ and $M=2I$, where $I$ is a $2\times 2$ identity matrix. Let
\begin{equation*}
\psi(x)=\left\{
	\begin{aligned}
	1 \quad\quad\quad\quad\quad  if\quad |x|\leq 1/3 \\
	cos[v(3|x|-1)\pi/2] \quad 1/3\leq|x|\leq 2/3\\
	0 \quad\quad\quad\quad\quad\quad otherwise\\
	\end{aligned}
	\right
	.
\end{equation*}
Here $v$ is a smooth and non-negative function satisfying \cite{TEN}
 \begin{equation*}
v(x)=\left\{
	\begin{array}{cl}
	0& \quad\quad\quad\quad\quad  if\quad x\leq 0 \\
    x^{4}(35-84x+70x^{2}-20x^{3})& \quad\quad\quad\quad\quad if\quad 0 <x< 1 \\
	1& \quad\quad\quad\quad\quad  if\quad x\geq 1\\
		\end{array}
	\right
	.
\end{equation*}
and $$v(x)+v(1-x)=1.$$
Let $\widehat{\tilde{\phi}}(w_{1},w_{2})=\psi(w_{1})\psi(w_{2})$, then $\tilde{\phi}\in W_{0}(L^{1}(\mathbb{R}^{2}))$, $\hat{\Phi}_{0}(w)=\sum_{k\in \mathbb{Z}^{2}}\widehat{\tilde{\phi}}(w+k)\neq 0$ and that $\{\tilde{\phi}(\cdot-k):k\in\mathbb{Z}^{2}\}$ forms a Riesz basis for $V(\tilde{\phi})$\cite[p137]{TEN}.

 \begin{figure}[h]
  \centering
  }
  \includegraphics[width=15cm,height=8cm]{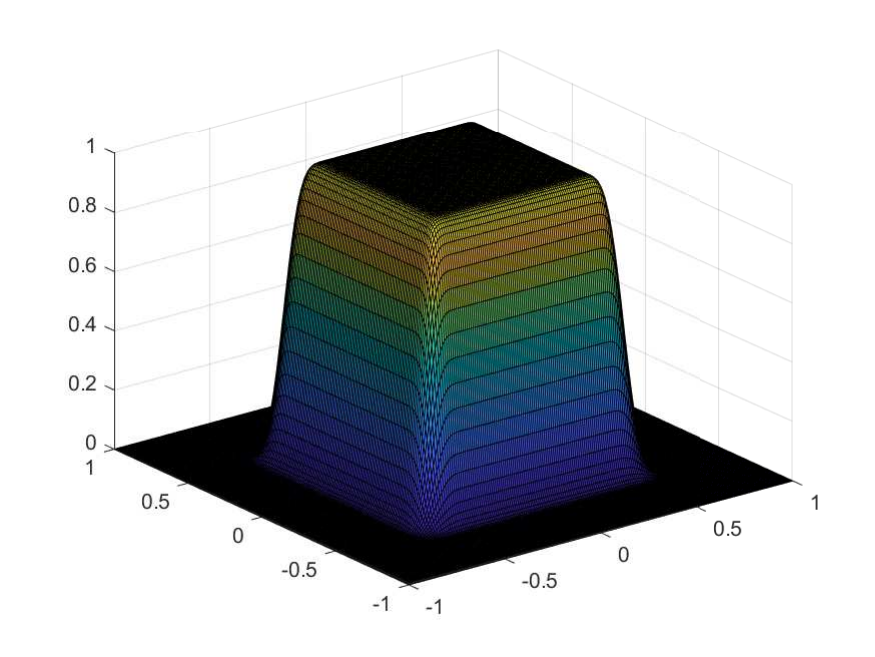}
  {
\caption{$\widehat{\widetilde{\phi}}$}
\end{figure}

Let sampled signal $f(\cdot)=\sum_{k\in\mathbb{Z}^{2}}c(k)\tilde{\phi}(\cdot-k)\in V(\tilde{\phi})$, where
\begin{equation*}
c(k)=\left\{
	\begin{array}{cl}
	1& \quad  (k_{1},k_{2})=(1,0) \\
	2& \quad  (k_{1},k_{2})=(0,1)\\
    0& \quad  otherwise\\
		\end{array}
	\right
.
\end{equation*}
Therefore, we can get $$f(w)=1\cdot\tilde{\phi}\left((w_{1},w_{2})-(1,0)\right)+2\cdot\tilde{\phi}\left((w_{1},w_{2})-(0,1)\right)=\tilde{\phi}(w_{1}-1,w_{2})+2\tilde{\phi}(w_{1},w_{2}-1).$$

\begin{figure}[htbp]
  \centering
  \begin{minipage}{0.49\linewidth}
  \centering
  \includegraphics[width=1.1\linewidth]{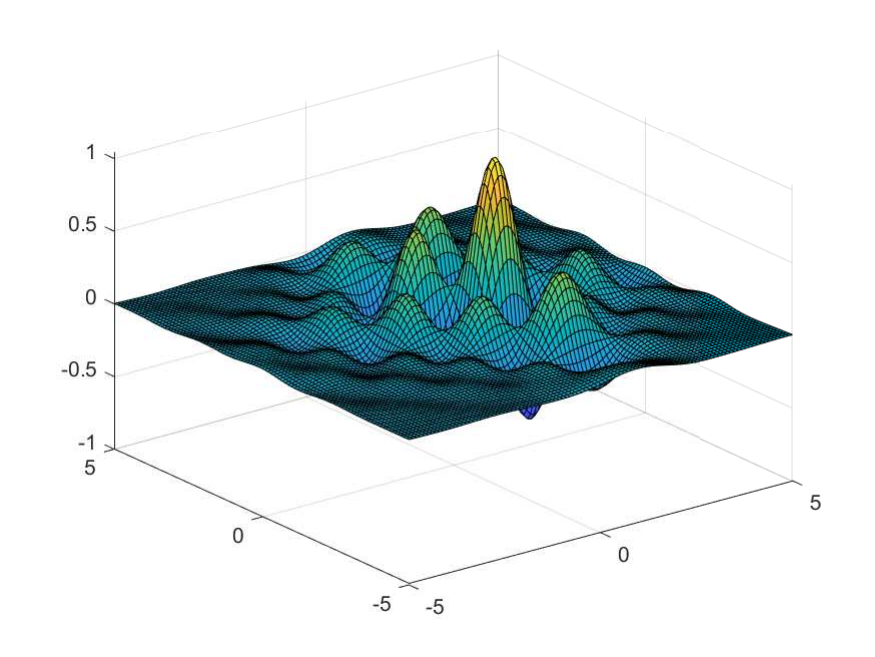}
  \caption{The real part of $f$.}
  \end{minipage}
  \begin{minipage}{0.49\linewidth}
  \centering
  \includegraphics[width=1.1\linewidth]{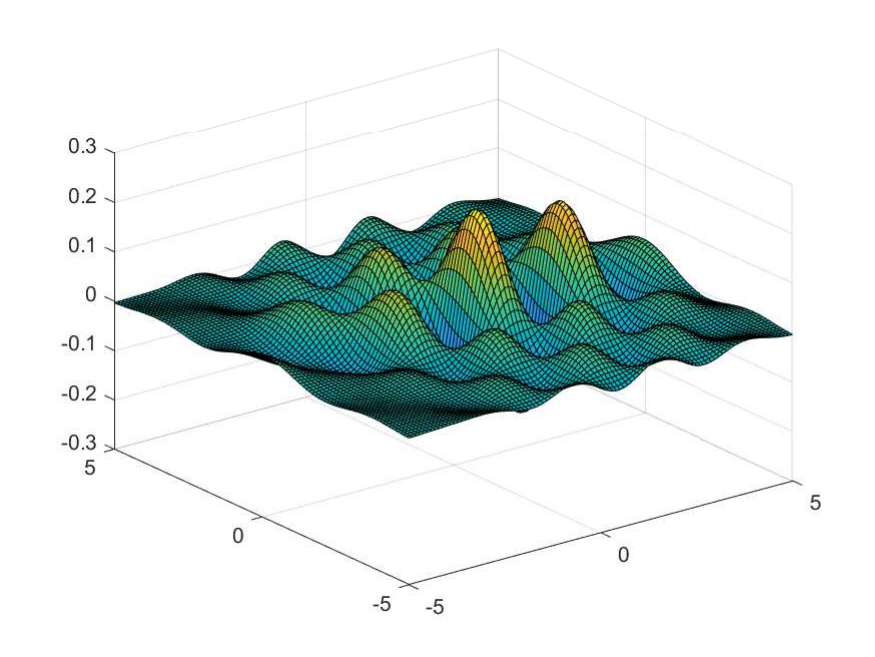}
  \caption{The imaginary part of $f$.}
  \end{minipage}
\end{figure}

\begin{figure}[htbp]
  \centering
  \begin{minipage}{0.49\linewidth}
  \centering
  \includegraphics[width=1.1\linewidth]{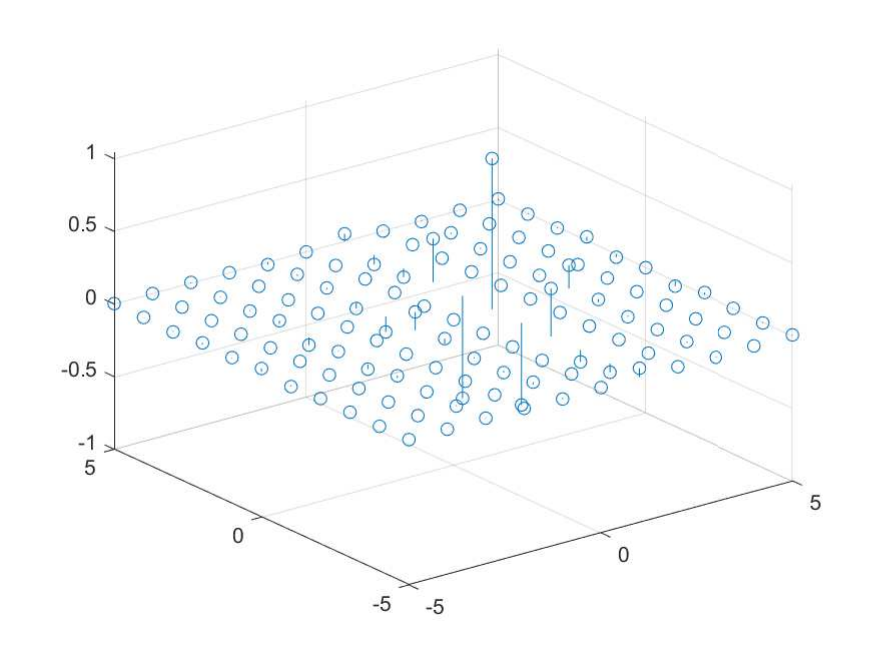}
  \caption{The real part of the sampling points of $f$.}
  \end{minipage}
  \begin{minipage}{0.49\linewidth}
  \centering
  \includegraphics[width=1.1\linewidth]{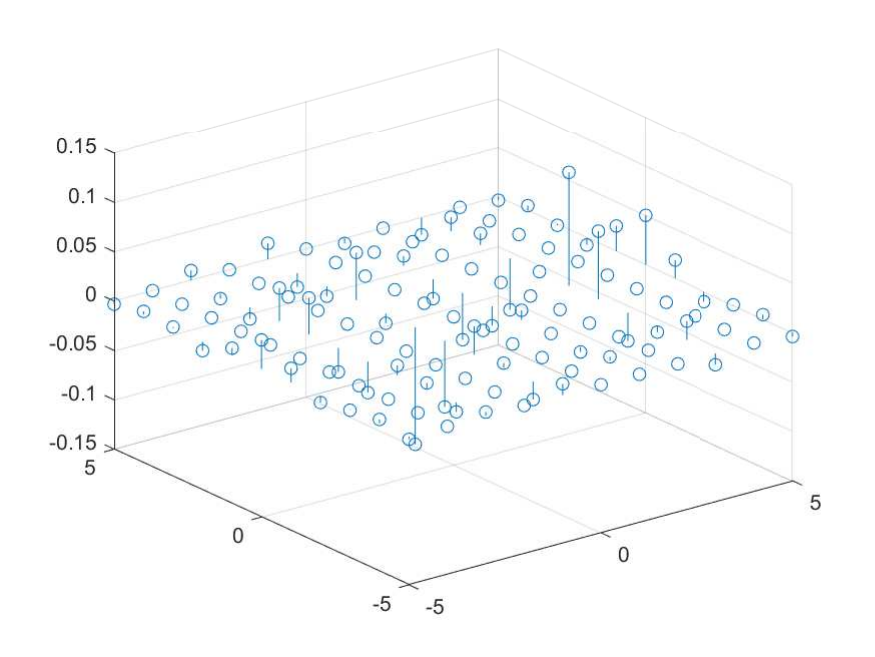}
  \caption{The imaginary part of the sampling points of $f$.}
  \end{minipage}
\end{figure}

Let
\begin{equation}
\phi(w)=\tilde{\phi}(w)e^{-i\pi [w^{T}B^{-1}Aw+2(B^{-1}P)^{T}w]}
\end{equation}
\begin{equation}
\tilde{f}(w)=f(w)e^{i\pi [w^{T}B^{-1}wt+2(B^{-1}P)^{T}w]}
\end{equation}

Meanwhile, we give the relation between the SAFT and the Fourier translation
\begin{equation}
\begin{aligned}\label{relation}
(S_\mathcal{M}f)(w)=&\frac{1}{\sqrt{|\mathrm{det}(B)}|}\int_{\mathbb{R}^{n}}f(t)e^{i\pi(t^{T}B^{-1}At+w^{T}DB^{-1}w-2t^{T}B^{-1}w)}\\
&\times e^{2i\pi(B^{-1}P)^{T}t+2i\pi(Q^{T}-P^{T}DB^{-1})w}dt\\
=&\frac{1}{\sqrt{|\mathrm{det}(B)}|}e^{i\pi[w^{T}DB^{-1}w+2(Q^{T}-P^{T}DB^{-1})w]}\widehat{\tilde{f}}(B^{-1}w).
\end{aligned}
\end{equation}

According to (\ref{qiuhe}) and (\ref{relation}), we can get
\begin{equation}
\begin{aligned}
(S_\mathcal{M}\Phi_{0})(w)&=\sum_{n\in\mathbb{Z}^{2}}\overline{\eta}^{2}_\mathcal{M}(w+n)(S_\mathcal{M}\phi)(w+n)\\
&=\sum_{n\in\mathbb{Z}^{2}}\frac{1}{\sqrt{|\mathrm{det}(B)}|}e^{-i\pi[w^{T}DB^{-1}w+2(Q^{T}-P^{T}DB^{-1})w]}\widehat{\tilde{\phi}}(B^{-1}(w+n)),
\end{aligned}
\end{equation}
where $\phi_{j}=a^{j}*\phi$.

\begin{figure}[htbp]
  \centering
  \begin{minipage}{0.49\linewidth}
  \centering
  \includegraphics[width=1.1\linewidth]{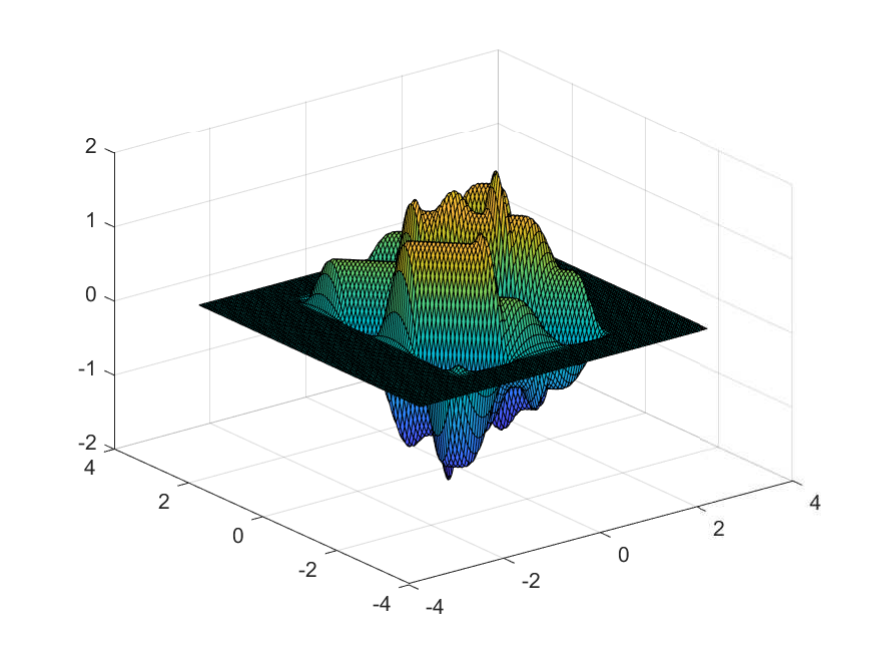}
  \caption{The real part of $(S_\mathcal{M}\Phi_{0})(w)$.}
  \end{minipage}
  \begin{minipage}{0.49\linewidth}
  \centering
  \includegraphics[width=1.1\linewidth]{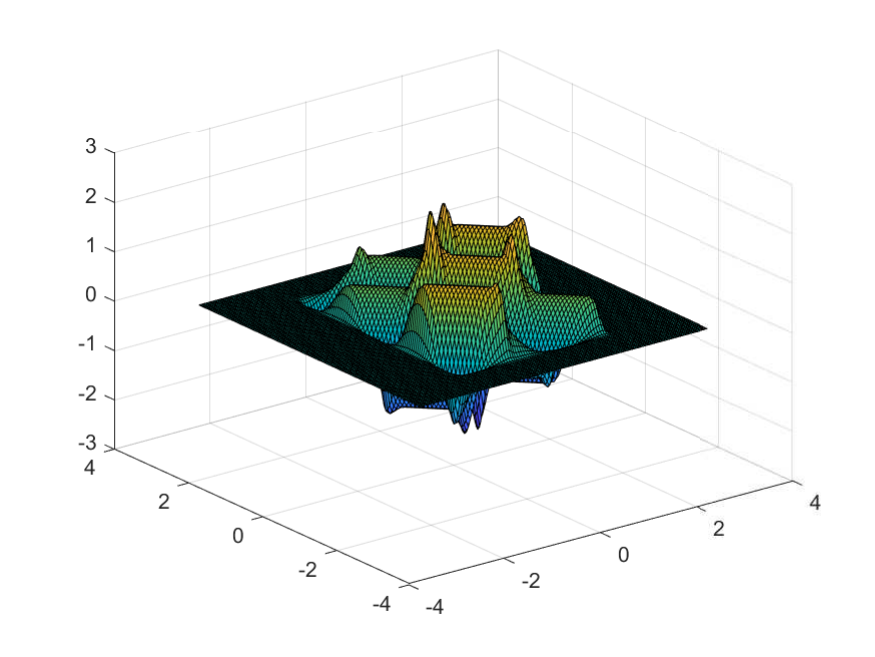}
  \caption{The imaginary part of $(S_\mathcal{M}\Phi_{0})(w)$.}
  \end{minipage}
\end{figure}

Let $\hat{a}(w)=\hat{b}(Bw)\chi_{E}(w)$, where $\hat{b}(w_{1},w_{2})=c_{1}e^{2\pi i(w_{1}+w_{2})}+c_{2}e^{2\pi i(w_{1}+2w_{2})}$.

\begin{figure}[htbp]
  \centering
  \begin{minipage}{0.49\linewidth}
  \centering
  \includegraphics[width=1.1\linewidth]{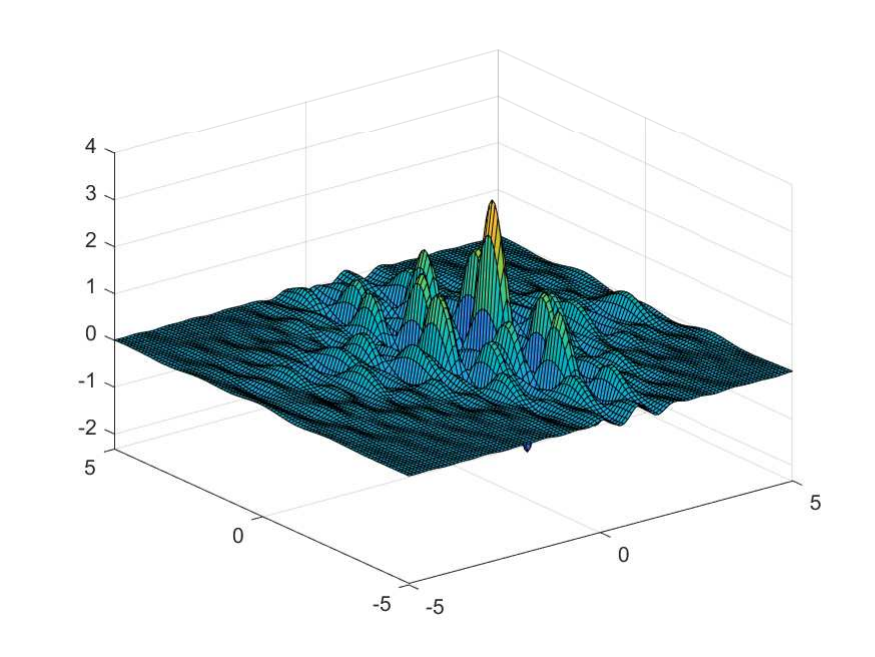}
  \caption{tupianThe real part of $a\ast\phi$.}
  \end{minipage}
  \begin{minipage}{0.49\linewidth}
  \centering
  \includegraphics[width=1.1\linewidth]{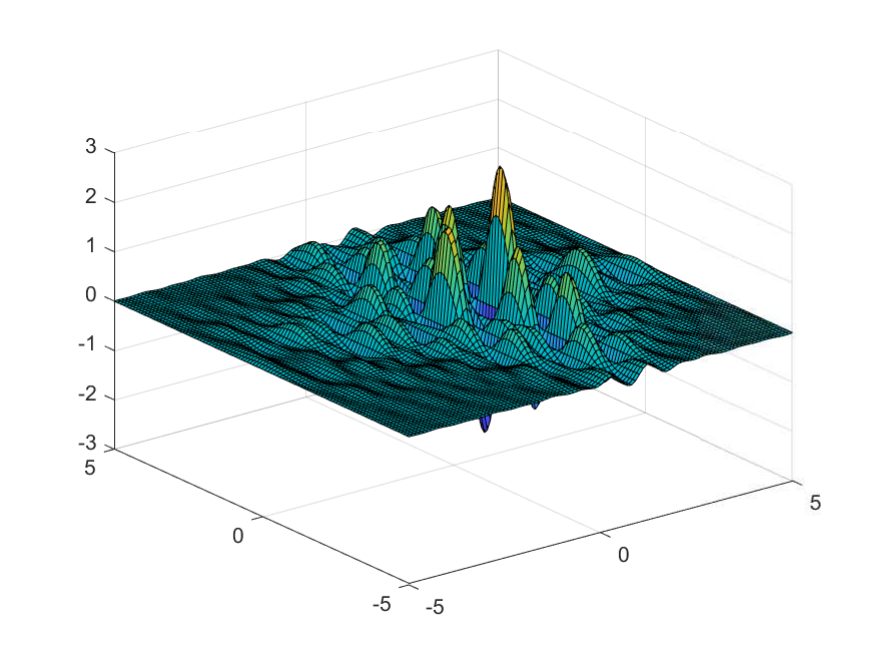}
  \caption{The imaginary part of $a\ast\phi$.}
  \end{minipage}
\end{figure}

\begin{figure}[htbp]
  \centering
  \begin{minipage}{0.49\linewidth}
  \centering
  \includegraphics[width=1.1\linewidth]{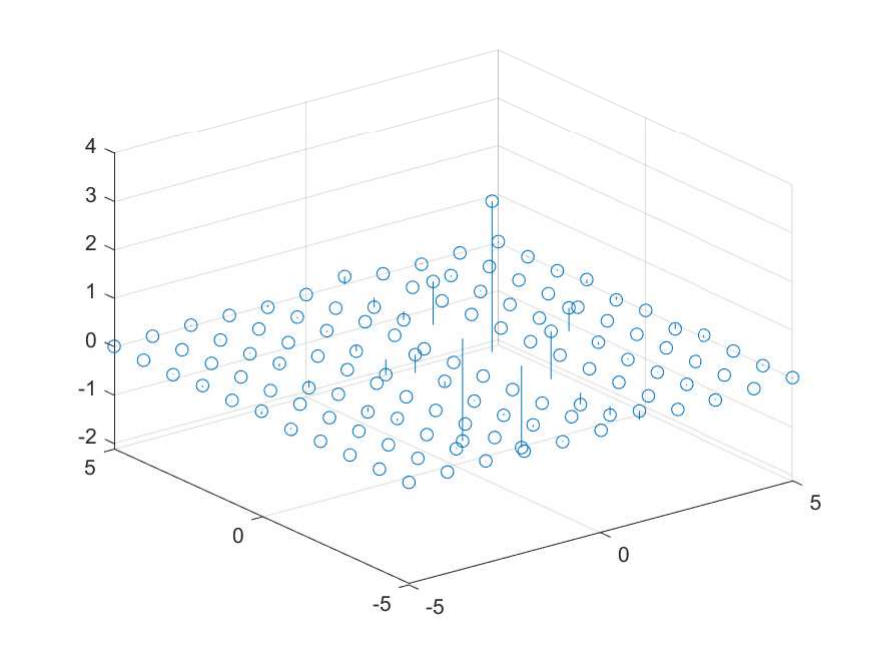}
  \caption{The real part of the sampling points of $a\ast\phi$.}
  \end{minipage}
  \begin{minipage}{0.49\linewidth}
  \centering
  \includegraphics[width=1.1\linewidth]{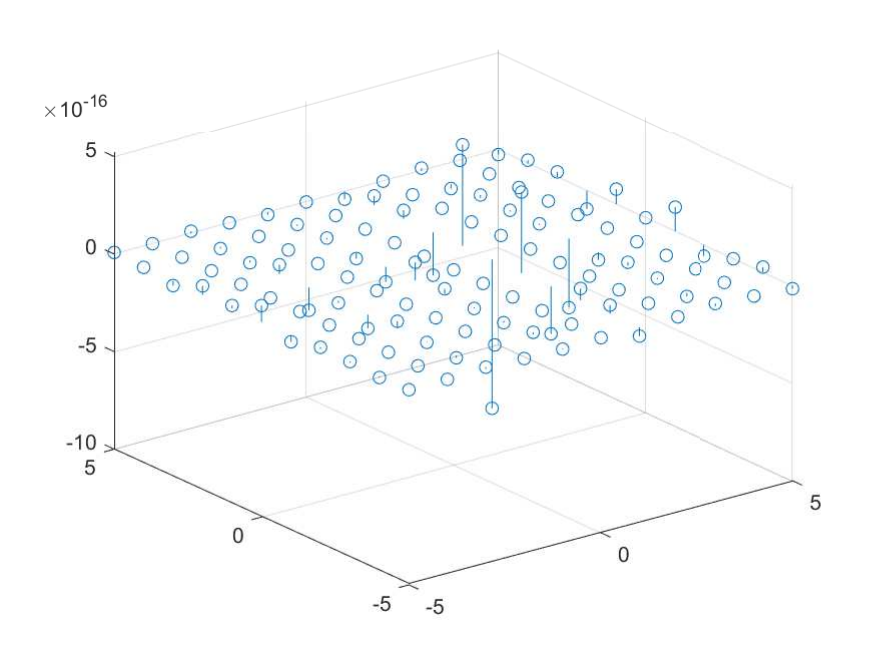}
  \caption{The imaginary part of the sampling points of $a\ast\phi$.}
  \end{minipage}
\end{figure}

Since for any $(k_{1},k_{2})^{T}\in\mathbb{Z}^{2}$ and $E=[-2/3,2/3]^{2}$, we have
\begin{eqnarray*}
&&\hat{a}(w_{1}+k_{1},w_{2}+k_{2})\widehat{\tilde{\phi}}(w_{1}+k_{1},w_{2}+k_{2})\\
&&\quad=\hat{b}(w_{1}+k_{1},w_{2}+k_{2})\chi_{E}(w_{1}+k_{1},w_{2}+k_{2})\widehat{\tilde{\phi}}(w_{1}+k_{1},w_{2}+k_{2})\\
&&\quad=\hat{b}(w_{1},w_{2})\chi_{E}(w_{1}+k_{1},w_{2}+k_{2})\widehat{\tilde{\phi}}(w_{1}+k_{1},w_{2}+k_{2})\\
&&\quad=\hat{b}(Bw)\widehat{\tilde{\phi}}(w_{1}+k_{1},w_{2}+k_{2}).
\end{eqnarray*}
That is to say,
\begin{equation}
\hat{a}(B^{-1}(w+n))\widehat{\tilde{\phi}}(B^{-1}(w+n))=\hat{b}(w)\widehat{\tilde{\phi}}(B^{-1}(w+n)).
\end{equation}
Then
\begin{equation}
\begin{aligned}
(S_\mathcal{M}\Phi_{1})(w_{1},w_{2})
&=\sum_{n\in\mathbb{Z}^{2}}\frac{1}{\sqrt{|\mathrm{det}(B)}|}e^{-i\pi[w^{T}DB^{-1}w+2(Q^{T}-P^{T}DB^{-1})w]}\widehat{\tilde{\phi}}(B^{-1}(w+n))\\
&=\sum_{n\in\mathbb{Z}^{2}}\frac{1}{\sqrt{|\mathrm{det}(B)}|}e^{-i\pi[w^{T}DB^{-1}w+2(Q^{T}-P^{T}DB^{-1})w]}(a*\tilde{\phi})^{\widehat{•}}(B^{-1}(w+n))\\
&=\sum_{n\in\mathbb{Z}^{2}}\frac{1}{\sqrt{|\mathrm{det}(B)}|}e^{-i\pi[w^{T}DB^{-1}w+2(Q^{T}-P^{T}DB^{-1})w]}\hat{a}(B^{-1}(w+n))\widehat{\tilde{\phi}}(B^{-1}(w+n))\\
&=\sum_{n\in\mathbb{Z}^{2}}\frac{1}{\sqrt{|\mathrm{det}(B)}|}e^{-i\pi[w^{T}DB^{-1}w+2(Q^{T}-P^{T}DB^{-1})w]}\hat{b}(w)\widehat{\tilde{\phi}}(B^{-1}(w+n))\\
&=\hat{b}(w)\sum_{n\in\mathbb{Z}^{2}}\frac{1}{\sqrt{|\mathrm{det}(B)}|}e^{-i\pi[w^{T}DB^{-1}w+2(Q^{T}-P^{T}DB^{-1})w]}\widehat{\tilde{\phi}}(B^{-1}(w+n))\\
&=\hat{b}(w)\sum_{n\in\mathbb{Z}^{2}}e^{-2i\pi[w^{T}DB^{-1}w+2(Q^{T}-P^{T}DB^{-1})w]}(S_\mathcal{M})\phi(w+n)\\
&=\hat{b}(w_{1},w_{2})(S_\mathcal{M}\Phi_{0})(w).
\end{aligned}
\end{equation}
So
\begin{equation}
\begin{aligned}
(S_\mathcal{M}\Phi_{1})(M^{-1}(w+\gamma_{k}))
=\hat{b}(M^{-1}(w+\gamma_{k}))(S_\mathcal{M}\Phi_{0})(M^{-1}(w+\gamma_{k}))
\end{aligned}
\end{equation}
Thus
\begin{eqnarray*}
\mathcal{D}(w)=\mathcal{E}(w)\mbox{diag}\left((S_\mathcal{M}\Phi_{0})\left(M^{-1}w\right),(S_\mathcal{M}\Phi_{0})\left(M^{-1}(w+\gamma_{1})\right)\right),
\end{eqnarray*}
where
\[
\mathcal{E}(w)=\left(
  \begin{array}{cccc}
   1& 1\\
    \widehat{b}\left(M^{-1}w\right)&\widehat{b}\left(M^{-1}(w+\gamma_{1})\right)\\
  \end{array}\right).
\]
Therefore $$\mbox{det}\,\mathcal{D}(w)=(S_\mathcal{M}\Phi_{0}))\left(M^{-1}w\right)(S_\mathcal{M}\Phi_{0})\left(M^{-1}(w+\gamma_{1})\right).$$
By $\mbox{det}\,\mathcal{E}(w)\neq0$. Hence $\mbox{det}\,\mathcal{D}(w)\neq0$. Thus, the $\hat{a}(w)$ which we chose makes $\mathcal{D}(w)$ satisfy the conditions in Theorem 4.6, that is, $\mathrm{det}\,\mathcal{D}(w)\neq0$.

\section{Conclusion}

In this paper, we investigate the multivariate dynamical sampling in the shift-invariant spaces
associated with the nD-SAFT. More specifically, we obtain
or a
function in a shift-invariant space $V (\phi)$ that is generated by $\phi\in l^{2}(\mathbb{R}^{n})$
can be stably recovered from its dynamical sampling measurements associated with the  nD-SAFT.  In the end, we give an example to verify the
effectiveness of our approach. Our results extend the original ones in the FT or LCT  domain.

\end{document}